\documentclass[11pt, twoside]{amsart}

\usepackage[T1]{fontenc}
\usepackage{amssymb,amsmath,amstext}

\newtheorem{theor}{Theorem}
\newtheorem{lem}[theor]{Lemma}

\newtheorem{exam}[theor]{Example}

\newtheorem{rem}[theor]{Remark}

\newcommand{\es}{\emptyset}

\newcommand{\mcA}{\mathcal{A}}
\newcommand{\mcB}{\mathcal{B}}

\newcommand{\mcF}{\mathcal{F}}
\newcommand{\mcG}{\mathcal{G}}

\newcommand{\mcL}{\mathcal{L}}
\newcommand{\mcM}{\mathcal{M}}
\newcommand{\mcN}{\mathcal{N}}

\newcommand{\mcP}{\mathcal{P}}

\newcommand{\mcT}{\mathcal{T}}

\newcommand{\mbC}{\mathbf{C}}

\newcommand{\mbG}{\mathbf{G}}

\newcommand{\mbP}{\mathbf{P}}

\newcommand{\mbS}{\mathbf{S}}

\newcommand{\mbbZ}{\mathbb{Z}}

\title[On compactness of logics]
{On compactness of logics that can express properties of symmetry or connectivity}

\author{Vera Koponen and Tapani Hyttinen}

\address{Vera Koponen, Department of Mathematics, Uppsala University, P.O. Box 480,
SE-75106 Uppsala, Sweden.}

\email{vera.koponen@math.uu.se}

\address{Tapani Hyttinen, Department of Mathematics, University of Helsinki, P.O. Box 68,
FI-00014 Helsinki, Finland}

\email{tapani.hyttinen@helsinki.fi} 

\begin{document}

\maketitle

\section{Introduction}

\noindent
A condition, in two variants, is given such that if a property $\mbP$ satisfies this condition,
then every logic which is at least as strong as first-order logic and can express $\mbP$
fails to have the compactness property.
The result is used to prove that for a number of natural properties $\mbP$ speaking about
automorphism groups or connectivity, 
every logic which is at least as strong as first-order logic and can express $\mbP$
fails to have the compactness property.
The basic idea underlying the results and examples presented here is that
it is possible to construct a countable first-order theory $\Gamma$ such that {\em every} model
of $\Gamma$ has a very rich automorphism group, but every finite $\Delta \subset \Gamma$ has
a model which is rigid. 
Although this can be proved by basic constructions, one can also derive the applications  
(Examples~\ref{1st example}--\ref{8th example}) 
in a uniform way by using the main theorems of this note, 
the proofs of which use results from random graph theory.
We conclude by showing that there is a logic that extends first-order logic,
has the compactness property and can
express the property ``the cardinality of the automorphism group is at most $2^{\aleph_0}$''.
What is common for all properties $\mbP$ of automorphism groups in Examples~\ref{1st example}--\ref{7th example}
is that $\mbP$ can be true as well as false for automorphism groups of cardinality at most $2^{\aleph_0}$.

Familiarity with first-order logic and basic model theory is assumed
(see \cite{EF, Rot} for example), as well as of basic group theory. 
For background about abstract logics see \cite{BF, Van}.
By $\sigma$ we denote a {\em countable} relational {\em signature}, also called {\em vocabulary}. 
In other words, $\sigma$ is a countable set
of relation symbols and with each relation symbol an arity is associated. 
Calligraphic letters $\mcA, \mcB, \mcM$, $\mcN$ etc. are used to denote structures and the corresponding
noncalligraphic letters $A$, $B$, $M$, $N$ etc. denote their universes.
Let $\mbS_\sigma$ denote the class of all $\sigma$-structures.

In this note we follow \cite{Van} and define an {\em abstract logic}, often just called a {\em logic}, 
over $\sigma$ to be
a pair $(L(\sigma), \models_{L(\sigma)})$, where $L(\sigma)$ is a set of objects, 
some of which are called sentences,
and $\models_{L(\sigma)} \subseteq \mbS_\sigma \times L(\sigma)$ is a relation such that
if $\varphi \in L(\sigma)$, $\mcM, \mcN \in \mbS_\sigma$ are isomorphic ($\mcM \cong \mcN$) and 
$\mcM \models_{L(\sigma)} \varphi$, then $\mcN \models_{L(\sigma)} \varphi$.
The set of first-order formulas that can be built from the vocabulary $\sigma$ is denoted $FO(\sigma)$ and
if $\models_{FO(\sigma)}$ is the usual satisfaction relation between $\sigma$-structures and first-order sentences,
then $(FO(\sigma), \models_{FO(\sigma)})$ is clearly an abstract logic.
From now on let $(L(\sigma), \models_{L(\sigma)})$ be an abstract logic over $\sigma$.
If $\Gamma \subseteq L(\sigma)$ and $\mcM \in \mbS_\sigma$ then $\mcM \models \Gamma$ means that
$\mcM \models_{L(\sigma)} \varphi$ for every $\varphi \in \Gamma$. 
The notation $\Gamma \models_{L(\sigma)} \varphi$
means that for every $\mcM \in \mbS_\sigma$, if $\mcM \models_{L(\sigma)} \Gamma$ 
then $\mcM \models_{L(\sigma)} \varphi$.
By a {\em property} $\mbP$ of $\sigma$-structures we mean a class $\mbP \subseteq \mbS_\sigma$ such that $\mbP$ is
closed under isomorphism. 
Let $\mbP$ be a property of $\sigma$-structures.
We say that {\em $\mbP$ can be expressed by $L(\sigma)$} if there is a sentence $\varphi \in L(\sigma)$ such that
for every $\mcM \in \mbS_\sigma$, $\mcM \models_{L(\sigma)} \varphi$ if and only if $\mcM \in \mbP$.
We say that {\em $L(\sigma)$ is at least as strong as $FO(\sigma)$}, in symbols $L(\sigma) \geq FO(\sigma)$,
if every property $\mbP \subseteq \mbS_\sigma$ that can be expressed by $FO(\sigma)$ can be expressed by $L(\sigma)$. 
Note that (for every property $\mbP \subseteq \mbS_\sigma$) the logic obtained by extending $FO(\sigma)$
with a generalised quantifier corresponding to $\mbP$ is at least as strong as $FO(\sigma)$ and can express $\mbP$.
If $L(\sigma) \geq FO(\sigma)$ and there is a property $\mbP \subseteq \mbS_\sigma$ such that $\mbP$ can be
expressed by $L(\sigma)$ but not by $FO(\sigma)$, then we write $L(\sigma) > FO(\sigma)$ and say that $L(\sigma)$ is stronger than $FO(\sigma)$.
We say that $L(\sigma)$ is {\em closed under negation} if for every sentence $\varphi \in L(\sigma)$
there is a sentence $\psi \in L(\sigma)$ such that for every $\mcM$, 
$\mcM \models_{L(\sigma)} \varphi$ if and only if $\mcM \not\models_{L(\sigma)} \psi$.

For any set $X$, $|X|$ denotes its cardinality. Let $\kappa$ be a cardinal.
We say that $L(\sigma)$ has the {\em $\kappa$-compactness property}
if the following holds:
\begin{itemize}
\item[] If $\Gamma \subseteq L(\sigma)$, $|\Gamma| \leq \kappa$ and for every finite $\Delta \subseteq \Gamma$,
there is $\mcM \in \mbS_\sigma$ such that $\mcM \models_{L(\sigma)} \Delta$, then there is 
$\mcM \in \mbS_\sigma$ such that $\mcM \models_{L(\sigma)} \Gamma$.
\end{itemize}
We say that $L(\sigma)$ has the {\em $\kappa$-Löwenheim-Skolem property} if the following holds:
\begin{itemize}
\item[] If $\Gamma \subseteq L(\sigma)$ and $\mcM \models_{L(\sigma)} \Gamma$, then there is
$\mcN \in \mbS_\sigma$ with universe $N$ of cardinality at most $\kappa$ such that $\mcN \models_{L(\sigma)} \Gamma$.
\end{itemize}
A basic fact about $FO(\sigma)$ is that it has both the  $\aleph_0$-compactness property and
the $\aleph_0$-Löwenheim-Skolem property.
Lindström's theorem \cite{Lin69}, as stated in \cite{Van}, says that if $L(\sigma) > FO(\sigma)$
then $L(\sigma)$ fails to have either the $\aleph_0$-compactness property or the $\aleph_0$-Löwenheim-Skolem
property. 
There exist abstract logics $L(\sigma) > FO(\sigma)$ with the $\aleph_0$-compactness property
as proved by Vaught \cite{Vau}, Keisler \cite{Kei}, Fuhrken \cite{Fuh}, Shelah \cite{She} and
Barwise, Kaufmann and Makkai \cite{BKM, BKM2}, for instance.

By a {\em proof system for $L(\sigma)$} we mean a triple $(\Pi, p, c)$ where
$\Pi$ is a set and $p$ and $c$ are functions with domain $\Pi$ such that
for every $\pi \in \Pi$, $c(\pi) \in L(\sigma)$ and $p(\pi)$ is a {\em finite} sequence of
members of $L(\sigma)$. For $\pi \in \Pi$ we call the members of $p(\pi)$ the {\em premisses} of
$\pi$ and $c(\pi)$ the {\em conclusion} of $\pi$. We allow the possibility that $p(\pi)$ is the empty sequence
(meaning that $\pi$ has no premisses).
A proof system $\Pi$ for $L(\sigma)$ is called
\begin{itemize}
\item[(i)] {\em sound} if every $\pi \in \Pi$ satisfies that if $p(\pi) = (\varphi_1, \ldots, \varphi_k)$,
$c(\pi) = \psi$ and $\varphi_1, \ldots, \varphi_k, \psi$ are sentences, then
$\{\varphi_1, \ldots, \varphi_k\} \models \psi$.
\item[(ii)] {\em complete} if whenever $\Gamma \subseteq L(\sigma)$ is a set of sentences,
$\psi \in L(\sigma)$ is a sentence and $\Gamma \models_{L(\sigma)} \psi$, then there is
$\pi \in \Pi$ such that $c(\pi) = \psi$ and every member of $p(\pi)$ belongs to $\Gamma$.
\end{itemize}

\noindent
The next section contains the main results, Theorems~\ref{main theorem} 
and~\ref{simplified form of main theorem}. In Section~\ref{examples},
Theorems~\ref{main theorem} 
and~\ref{simplified form of main theorem}
are applied to give a number of examples of properties $\mbP \subseteq \mbS_\sigma$ 
such that every logic $\mcL(\sigma) \geq FO(\sigma)$
which can express $\mbP$ is not $\aleph_0$-compact and (if it is closed under negation) does not have
a sound and complete proof system. It follows that such $\mbP$ cannot be expressed by $FO(\sigma)$.
In Section~\ref{proof of main theorem} Theorem~\ref{main theorem} is proved. The proof of 
Theorem~\ref{simplified form of main theorem} is a simpler variation of the 
proof of Theorem~\ref{main theorem}.
Section~\ref{Proofs of the applications without random graph theory} explains how 
the conclusions of Examples~\ref{1st example}--\ref{8th example} can be derived without
the use of random graph theory.
The final section proves that there is a logic which extends first-order logic, is $\aleph_0$-compact
and can express the property ``the cardinality of the automorphism group is at most $2^{\aleph_0}$''.

\section{Results}\label{results}

\noindent
For the rest of this section we fix a countable relational signature $\sigma$ and 
an abstract logic $(L(\sigma), \models_{L(\sigma)})$
which we abbreviate by $L(\sigma)$.
We will use notions of a graph theoretic flavour.
Let $\mcM$ be a $\sigma$-structure.
For $a_1, a_2 \in M$ we say that $a_1$ and $a_2$ are {\em adjacent (in $\mcM$)},
denoted $a_1 \sim_\mcM a_2$, if $a_1 \neq a_2$ and there are a relation symbol $R \in \sigma$,
$(b_1, \ldots b_r) \in R^\mcM$ and $i, j \in \{1, \ldots, r\}$ such that $a_1 = b_i$ and $a_2 = b_j$.
A sequence of distinct $a_1, \ldots, a_k \in \mcM$ (where $k$ is a positive integer) is called a 
{\em path (in $\mcM$)} if $a_i \sim_\mcM a_{i+1}$ for all $i \in \{1, \ldots, k-1\}$.
A sequence $a_1, \ldots, a_k \in M$ such that $k \geq 3$, $a_1 = a_k$, $a_{k-1} \sim_\mcM a_k$
and $a_1, \ldots, a_{k-1}$ is a path is called a 
{\em cycle of length $k$ (in $\mcM$)} or a {\em $k$-cycle}. 
We say that $\mcM$ is {\em connected} if for all distinct $a_1, a_2 \in M$ there is a path
$b_1, \ldots, b_k \in M$ such that $a_1 = b_1$ and $a_2 = b_k$.
The {\em distance} between $a$ and $b$ is 0 if $a = b$ and otherwise it is the the minimal
integer $k$ such that there is a path starting with $a$ and ending with $b$ that contains $k+1$ elements;
if no such $k$ exists then we say that the distance between $a$ and $b$ is $\infty$.
For $a \in M$ the {\em degree of $a$ (in $\mcM$)} is the cardinality of
\[ \{b \in M : b \sim_\mcM a\}. \]
Let $R \in \sigma$ be a relation symbol of arity $r \geq 2$.
By an {\em $R$-graph} we mean a $\sigma$-structure $\mcM$ such that
$Q^\mcM = \es$ for every $Q \in \sigma \setminus \{R\}$, and 
for all $(a_1, \ldots, a_r) \in M^r$,
\[
(a_1, a_2, \ldots, a_r) \in R^\mcM \ \Longrightarrow \ a_1 \neq a_2, \ 
a_2 = \ldots = a_r \ \text{ and } 
(a_2, \underbrace{a_1, \ldots, a_1}_{r-1 \text{ times}}) \in R^\mcM.
\]
Note that if $\mcM$ is an $R$-graph and
\[ E \ = \ \Big\{ \{a, b\} \subseteq M : (a,b, \ldots, b) \in R^\mcM \Big\} \]
then $(M, E)$ is an undirected graph without loops in the usual sense of graph theory 
(see \cite{Die} for example).
Conversely, if $(M, E)$ is an undirected graph without loops and 
\begin{align*}
&R^\mcM \ = \ 
\Big\{ (a_1, a_2,  \ldots, a_r) \in M^r : a_2 = \ldots = a_r \text{ and } \{a_1, a_2\} \in E \Big\} 
\ \text{ and} \\
&Q^\mcM \ = \ \es \ \text{ for every } Q \in \sigma \setminus \{R\},
\end{align*}
then $\mcM$ is an $R$-graph.
An $R$-graph is called an {\em $R$-tree} if it is connected and has no cycle.

Suppose that $I$ is a set and, for every $i \in I$, that $\mcM_i$ is a $\sigma$-structure.
Then $\bigcup_{i \in I} \mcM_i$ denotes the $\sigma$-structure $\mcM$ with universe 
$M = \bigcup_{i \in I} M_i$ and such that, for every $Q \in \sigma$,
$Q^\mcM = \bigcup_{i \in I} Q^{\mcM_i}$.
By $\mcM_1 \cup \mcM_2$ we mean $\bigcup_{i \in \{1,2\}} \mcM_i$.
Let $R \in \sigma$ have arity at least 2. Then $\mcT_R$ denotes an infinite $R$-tree
such that exactly one element of its universe has degree $4$ and all other elements of its universe have degree 5.
$\mcT'_R$ denotes the infinite $R$-tree such that all elements of its universe have degree 5.

For every $\sigma$-structure $\mcM$,
let $Aut(\mcM)$ denote the group of automorphisms of $\mcM$. 
A $\sigma$-structure $\mcM$ is called {\em rigid}
if $Aut(\mcM)$ contains only one element (the identity function).

If $Q \in \sigma$ and $\mcM$ is a $\sigma$-structure such that $Q^\mcM = \es$,
then $\mcM_Q$ denotes the $\sigma$-structure with universe $M_Q = M$ such that 
\[ Q^{\mcM_Q} \ = \  \{(a_1, \ldots, a_q) \in M^q : a_1 = \ldots = a_q\}, \]
where $q$ is the arity of $Q$, and $R^{\mcM_Q} = R^\mcM$ for all $R \in \sigma \setminus \{Q\}$.
Recall that $\mbS_\sigma$ denotes the class of all $\sigma$-structures.
For $\mcM, \mcN \in \mbS_\sigma$, the notation $\mcM \equiv \mcN$ means, as usual, that
$\mcM$ and $\mcN$ are elementarily equivalent, i.e., satisfy exactly the same $FO(\sigma)$-sentences.

\begin{theor}\label{main theorem}
Suppose that $Q, R \in \sigma$ are different relation symbols and that $R$ has arity at least 2.
Let $\mbP \subseteq \mbS_\sigma$ be a property (so it is closed under isomorphism)
such that there is $\mcM \in \mbS_\sigma$ such that $Q^\mcM = \es$ and the following conditions hold:
\begin{itemize}
\item[(1)] If $\mcN \in \mbS_\sigma$ is a finite, rigid, connected, $N \neq \es$,
$Q^\mcN = \es$ and $M \cap N = \es$, then $\mcM_Q  \cup  \mcN \ \in \ \mbP$.

\item[(2)] If $\mcM' \equiv \mcM$ and 
$I$ and $J$ are sets such that $I$ is infinite and for all $i \in I$ and all $j \in J$,
$\mcT_i \cong \mcT_R$, $\mcT'_j \cong \mcT'_R$ 
and the universes of any two structures from $\{\mcM'\} \cup \{\mcT_i : i \in I\} \cup \{\mcT'_j : j \in J\}$ 
are disjoint, then
\[ \mcM'_Q \ \cup \ \bigcup_{i \in I} \mcT_i \ \cup \ \bigcup_{j \in J} \mcT'_j \ \notin \ \mbP. \]
\end{itemize}
Then the following hold:
\begin{itemize}
\item[(i)] $\mbP$ cannot be expressed by $FO(\sigma)$.
\item[(ii)] If $L(\sigma) \geq FO(\sigma)$ and $\mbP$ can be expressed by $L(\sigma)$, 
then $L(\sigma)$ does not have the $\aleph_0$-compactness property.
\item[(iii)] If $L(\sigma) \geq FO(\sigma)$, $L(\sigma)$ is closed under negation and
$\mbP$ can be expressed by $L(\sigma)$, then there is no proof system for $L(\sigma)$
which is both sound and complete.
\end{itemize}
\end{theor}

\noindent
Below follows a variant of the above result which is applicable also in the case when the signature
contains only one relation symbol.

\begin{theor}\label{simplified form of main theorem}
Suppose that $R \in \sigma$ is a relation symbol with arity at least 2.
Let $\mbP \subseteq \mbS_\sigma$ be a property such that the following conditions hold:
\begin{itemize}
\item[(1)] If $\mcN \in \mbS_\sigma$ is finite, rigid, connected and $N \neq \es$,
then $\mcN \ \in \ \mbP$.

\item[(2)] If $I$ and $J$ are sets such that $I$ is infinite and for all $i \in I$ and all $j \in J$,
$\mcT_i \cong \mcT_R$, $\mcT'_j \cong \mcT'_R$ 
and the universes of any two structures from $\{\mcT_i : i \in I\} \cup \{\mcT'_j : j \in J\}$ 
are disjoint, then
\[ \bigcup_{i \in I} \mcT_i \ \cup \ \bigcup_{j \in J} \mcT'_j \ \notin \ \mbP. \]
\end{itemize}
Then~(i)--(iii) of Theorem~\ref{main theorem} hold.
\end{theor}

\noindent
The next section gives applications of the above theorems.
Theorem~\ref{main theorem} is proved in Section~\ref{proof of main theorem}.
The proof of Theorem~\ref{simplified form of main theorem} is a 
(simpler) modification of the proof of Theorem~\ref{main theorem}.

Let us first put the theorems above in some perspective.
It is easy to see, and well known, that if $L(\sigma) \geq FO(\sigma)$ is a logic
which can express the property 
\[ \mbP_f \ = \ \{\mcM \in \mbS_\sigma : M \text{ is finite} \}, \]
then $L(\sigma)$ is not $\aleph_0$-compact.
Clearly~(1) and~(2) of Theorem~\ref{simplified form of main theorem} hold, so the
theorem gives the expected result.
Now let
\[ \mbP_c \ = \ \{\mcM \in \mbS_\sigma : M \text{ is countable} \}.\]
As proved by Vaught \cite{Vau} (see also \cite{Fuh, Kei}), 
the extension of $FO(\sigma)$ by a generalised
quantifier that expresses ``there are uncountably many'' is $\aleph_0$-compact
and can express $\mbP_c$ (as this logic is closed under negation). 
Therefore~(1) or~(2) of Theorem~\ref{simplified form of main theorem}
must fail. Since~(1) clearly holds for $\mbP = \mbP_c$,~(2) must fail,
which is also easy to see directly, because if $I \cup J$ is countable then
the corresponding structure belongs to $\mbP_c$, while if 
$I \cup J$ is uncountable then the resulting structure does not belong to $\mbP_c$.

\section{Applications}\label{examples}

\noindent
In this section examples are given of properties $\mbP$ for which conditions~(1) and~(2)
of Theorem~\ref{main theorem} or Theorem~\ref{simplified form of main theorem} are satisfied.
Hence we get examples of properties $\mbP$ such that every logic $L(\sigma) \geq FO(\sigma)$
in which $\mbP$ can be expressed is not $\aleph_0$-compact and (if $L(\sigma)$ is closed 
under negation) does not have a sound and complete proof system.
Most of the properties considered are properties of the group of automorphisms of a structure
(where the group operation is composition). For a structure $\mcM$ the group of automorphisms,
or automorphism goup, is denoted $Aut(\mcM)$.
Automorphism groups, or symmetry groups as they are called in geometry, are of interest
since their properties are related to the properties of the underlying structures. 
In the words of P. de la Harpe \cite{Har}: 
``{\em ... symmetries and groups is one way of coping with the frustrations of life's limitations:
we like to recognize symmetries which allow us to grasp more than what we can see.}''
Moreover, the study of automorphism groups provides examples, concepts and methods for abstract group theory.
(See \cite{Eva, KM} for automorphism groups of first-order structures.)

The lemma below will simplify applications of Theorems~\ref{main theorem} 
and~\ref{simplified form of main theorem}.
The basic idea is that if $\mbP$ is a property that, for countable $\mcM \in \mbP$, restricts the
complexity of $Aut(\mcM)$, then~(1) and~(2) of either 
Theorem~\ref{main theorem} or Theorem~\ref{simplified form of main theorem} hold.

Suppose that the signature $\sigma$ has a relation symbol $R$ with arity at least 2,
so it makes sense to talk about the $R$-trees $\mcT_R$ and $\mcT'_R$ from the previous section.
For every nonempty set $X$ let $S(X)$ be the {\em symmetric group of $X$}, that is,
the group of all bijective functions (i.e. permutations) from $X$ to itself.
Let $\mcG$ be the automorphism group of the $\sigma$-structure
\[ \mcM'_Q \ \cup \ \bigcup_{i \in I} \mcT_i \ \cup \ \bigcup_{j \in J} \mcT'_j \]
from~(2) of Theorem~\ref{main theorem} or of the $\sigma$-structure
$\bigcup_{i \in I} \mcT_i \ \cup \ \bigcup_{j \in J} \mcT'_j$
from~(2) of Theorem~\ref{simplified form of main theorem}.
Since $\mcT_i \cong \mcT_j$ for all $i,j \in I$ it follows that $S(I)$ is isomorphic with a
subgroup of $\mcG$. Note that (by assumption in~(2)) $I$ is infinite.
In particular, it follows that the cardinality of $\mcG$ is at least $2^{\aleph_0}$ and
that $\mcG$ has a subgroup which is isomorphic with $S(\omega)$, where $\omega$ is the set of
finite cardinals, so $\mcG$ is a quite ``rich'' group. 
By (the proof of) Cayley's theorem, every countable group is isomorphic with a subgroup of $\mcG$.

\begin{lem}\label{application lemma}
Suppose that $\mbP \subseteq \mbS_\sigma$ is a property such that if 
$\mcM, \mcN \in \mbS_\sigma$ and $Aut(\mcM) \cong Aut(\mcN)$, then
$\mcM \in \mbP$ if and only if $\mcN \in \mbP$.
\begin{itemize}
\item[(a)] Suppose that $\mbP$ is such that 
if $Aut(\mcM)$ has a subgroup which is isomorphic with 
$S(\omega)$, then $\mcM \notin \mbP$. 
Then~(2) of Theorem~\ref{main theorem} and of Theorem~\ref{simplified form of main theorem}
is satisfied.

\item[(b)] If every rigid structure $\mcM \in \mbS_\sigma$ belongs to $\mbP$, 
then~(1) of Theorem~\ref{simplified form of main theorem}
is satisfied.

\item[(c)] Suppose that $\sigma$ contains a relation symbol $Q \neq R$.
If there is a structure $\mcM \in \mbS_\sigma \cap \mbP$ such that
$Q^\mcM = \es$ then~(1) of Theorem~\ref{main theorem} is satsified.
\end{itemize}
\end{lem}

\noindent
{\bf Proof.}
Suppose that $\mbP$ depends only on the isomorphism type of $Aut(\mcM)$.
Part~(a) follows from what was said in the paragraph before the lemma.
If every rigid $\mcM \in \mbS_\sigma$ belongs to $\mbP$, then clearly 
the implication stated by~(1) of Theorem~\ref{simplified form of main theorem} holds.
For~(c), suppose that $\mcM \in \mbS_\sigma \cap \mbP$ and $Q^\mcM = \es$.
Also suppose that $\mcN \in \mbS_\sigma$ is finite, rigid, connected, $N \neq \es$, $Q^\mcN = \es$
and $M \cap N = \es$.
Consider $Aut(\mcM_Q \cup \mcN)$. Note that $\mcM_Q \cup \mcN \models Q(a, \ldots, a)$ if $a \in M$
and $\mcM_Q \not\models Q(a, \ldots, a)$ if $a \in N$, so automorphisms of $\mcM_Q \cup \mcN$
cannot send elements from $N$ to $M$ or vice versa. Since $Aut(\mcM) = Aut(\mcM_Q)$
and $\mcN$ is rigid it follows that $Aut(\mcM_Q \cup \mcN) \cong Aut(\mcM)$ so 
$\mcM_Q \cup \mcN \in \mbP$. So~(c) is proved.
\hfill $\square$

\begin{exam}\label{1st example}{\rm
Suppose that $\mbC$ is a class of groups such that for every $\mcG \in \mbC$,
$S(\omega)$ is not a subgroup of $\mcG$. 
Let
\[ \mbP \ = \ \{\mcA \in \mbS_\sigma : Aut(\mcA) \cong \mcG \text{ for some } \mcG \in \mbC \}. \]
By Lemma~\ref{application lemma}~(1), $\mbP$ satisfies~(2) of Theorems~\ref{main theorem}
and~\ref{simplified form of main theorem}.
If the trivial group belongs to $\mbC$ or if there is $\mcG \in \mbC$ such that
$\mcG \cong Aut(\mcM)$ for some $\mcM \in \mbS_\sigma$ with $Q^\mcM = \es$,
then $\mbP$ satisfies~(1) of Theorem~\ref{main theorem} or~\ref{simplified form of main theorem}.
So the conclusions (i)--(iii) of either theorem holds.
}\end{exam}

\begin{exam}\label{2nd example}{\rm
Suppose that $\mbC$ is a class of groups such that every member of $\mbC$ has cardinality less than $2^{\aleph_0}$
and define $\mbP$ as in Example~\ref{1st example}.
Since $S(\omega)$ has cardinality $2^{\aleph_0}$ it follows that no member of $\mbC$ has a subgroup
which is isomorphic with $S(\omega)$. 
Hence the conclusions of Example~\ref{1st example} are applicable in this case.
}\end{exam}

\begin{exam}\label{3rd example}{\rm
Suppose that $\Omega$ is a set of cardinals all of which are smaller than $2^{\aleph_0}$.
Let 
\[\mbP \ = \ \{\mcA \in \mbS_\sigma : \text{ the cardinality of $Aut(\mcA)$ belongs to $\Omega$}\}.\]
If we let $\mbC$ be the class of groups whose cardinality belongs to $\Omega$, then
\[\mbP \ = \ \{\mcA \in \mbS_\sigma : Aut(\mcA) \cong \mcG \text{ for some } \mcG \in \mbC \}, \]
so the conclusions of Example~\ref{1st example} are applicable.
In particular the properties ``the automorphism group is countable'' and
``the automorphism group is finite'' satisfies~(1) and~(2) of Theorems~\ref{main theorem}
and~\ref{simplified form of main theorem}.
}\end{exam}

\begin{exam}\label{4th example}{\rm
Every finitely generated, or countably generated group, is countable.
In particular, every finitely presentable group is countable.
Since $S(\omega)$ is not countable and the trivial group is finitely generated and finitely presentable
it follows from Example~\ref{2nd example} that the properties
``the automorphism group is finitely generated'',
``the automorphism group is countably generated'' and ``the automorphism group is finitely presentable''
satisfy~(1) and~(2) of Theorem~\ref{simplified form of main theorem}.
The argument can easily be modified to give similar conclusions for properties like
``the autmorphism group is nontrivial and finitely generated'' or 
``the automorphism group is infinite and finitely generated''.
}\end{exam}

\begin{exam}\label{5th example}{\rm
A group is called {\em locally finite} if every finitely generated subgroup of it is finite.
Since the group $\mbbZ$ is isomorphic with a subgroup of $S(\omega)$ it follows that
if a group has a subgroup which is isomorphic with $S(\omega)$, then it is not locally finite.
Since the trivial group is locally finite it follows from Example~\ref{1st example}
that~(1) and~(2) of Theorem~\ref{simplified form of main theorem} hold.
As one can find (exercise) $\mcM \in \mbS_\sigma$ such that $Aut(\mcM)$ is
infinite, locally finite and $Q^\mcM = \es$,
it follows that the property ``the automorphism group is infinite and
locally finite'' satisfies~(1) and~(2) of Theorem~\ref{main theorem}.
}\end{exam}

\noindent
The remaining examples do not consider properties that only depend on the automorphism group
as an abstract group, so Lemma~\ref{application lemma} is not applicable.

\begin{exam}\label{6th example}{\rm
The {\em support} of a function $f$, denoted $s(f)$, is the set of elements $a$ of the domain of $f$
such that $f(a) \neq a$.
Suppose that $Q, R \in \sigma$ are different where the arity of $R$ is at least 2. 
Let $X \subseteq \omega$ be such that there exists $\mcM \in \mbS_\sigma$ such that $Q^\mcM = \es$
and $|s(f)| \in X$ for every $f \in Aut(\mcG)$ and note that $Aut(\mcM) = Aut(\mcM_Q)$.
Let
\[ \mbP \ = \ \{\mcA \in \mbS_\sigma : \text{$|s(f)| \in X$ for every $f \in Aut(\mcA)$}\}. \]
If $\mcN \in \mbS_\sigma$ is rigid and $M \cap N = \es$
then $\mcM_Q \cup \mcN \in \mbP$. 
Hence~(1) of Theorem~\ref{main theorem} is satisfied.
As observed in the beginning of this section, for every structure as in~(2) of
Theorem~\ref{main theorem}, its automorphism group has a subgroup which is isomorphic with $S(\omega)$
and therefore it has (infinitely many) automorphisms with infinite support.
Hence~(2) of the same theorem is satisfied.
One can easily adapt the argument so that show that~(1) and~(2) of Theorem~\ref{main theorem} 
or~\ref{simplified form of main theorem}
hold for variations of the property defined in this example.
}\end{exam}

\begin{exam}\label{7th example}{\rm

An {\em orbit} of a structure $\mcA$ is a set $O \subseteq A$ such that whenever
$a, b \in O$, then there is $f \in Aut(\mcA)$ such that $f(a) = b$ and if
$a \in O$, $b \in A$ and $f(a) = b$ for some $f \in Aut(\mcA)$, then $b \in O$.
Hence the set of orbits of $\mcA$ is a partition of $A$. 
Let $o(\mcA)$ denote the set of orbits of $\mcA$.
Suppose that $Q, R \in \sigma$ are different where the arity of $R$ is at least 2. 
Consider the structure in~(2) of Theorem~\ref{main theorem}:
\[ \mcM'_Q \ \cup \ \bigcup_{i \in I} \mcT_i \ \cup \ \bigcup_{j \in J} \mcT'_j. \]
We know nothing more about $\mcM'$ than that $\mcM' \equiv \mcM$, 
so even if we know some properties of $Aut(\mcM)$
we do not know what $Aut(\mcM')$ is like, and hence not what $Aut(\mcM'_Q)$ is like.
It is clear that $\mcT_i \not\cong \mcT'_j$ for all $i \in I$ and $j \in J$.
The substructure $\bigcup_{i \in I} \mcT_i$ has exactly $\aleph_0$ orbits, because 
two elements in this structure are in the same orbit if and only if the distance 
(in the obvious sense) to the unique element with degree 4 in their respective connected component is the same.
Moreover, every orbit of $\bigcup_{i \in I} \mcT_i$ is infinite.
If $J \neq \es$, then the substructure $\bigcup_{j \in J} \mcT'_j$ has exactly 1 orbit, because
every $\mcT'_j$ is an infinite $R$-tree in which every element has degree 5.
Hence, the set of all orbits of the entire structure has cardinality $\aleph_0 + |o(\mcM')|$.

Note that whenever $\mcA, \mcN \in \mbS_\sigma$, $Q^\mcA = Q^\mcN = \es$, 
$\mcN$ is finite, rigid, $N \neq \es$  and $A \cap N = \es$, then all orbits of $\mcN$ are singletons,
so the cardinality of the set of orbits of $|o(\mcA_Q \cup \mcN)| = |o(\mcA)| + |N|$. 

Using these observations one can show, by choosing appropriate $\mcM \in \mbS_\sigma$ in each case,
that each of the following properties, to mention some examples, satisfy~(1) and~(2) of Theorem~\ref{main theorem}:
\begin{align*}
&\{\mcA \in \mbS_\sigma : |o(\mcA)| < \aleph_0\},\\
&\{\mcA \in \mbS_\sigma : \mcA \text{ has at least $k$ finite orbits}\}, \ \text{ for $k < \aleph_0$ },\\
&\{\mcA \in \mbS_\sigma : \mcA \text{ has only orbits of cardinality $< \kappa$}\}, \ 
\text{ for $\kappa \leq \aleph_0$}, \text{ and} \\
&\{\mcA \in \mbS_\sigma : \mcA \text{ has at most $k$ orbits of cardinality $> 1$}\}
\ \text{ for $k < \aleph_0$}.
\end{align*}
}\end{exam}

\noindent
We conclude this section with an example which does not speak about automorphism groups.

\begin{exam}\label{8th example}{\rm
Suppose that $Q, R \in \sigma$ are different where the arity of $R$ is at least 2. 
Clearly the structure in~(2) of Theorem~\ref{main theorem} have infinitely many connected components.
For any positive integer $k$ we can choose $\mcM$ with $Q^\mcM = \es$ and exactly
$k$ connected components. For any connected (nonempty) $\mcN$ such that $M \cap N = \es$, 
$\mcM_Q \cup \mcN$ has exactly $k+1$ components. Hence for any integer $k > 1$ and appropriate $\mcM$,
the property
\[ \{ \mcA \in \mbS_\sigma : \mcA \text{ has at most $k$ components} \} \]
satisfies~(1) and~(2) of Theorem~\ref{main theorem}.
For $k = 1$ the above property satisfies~(1) and~(2) of Theorem~\ref{simplified form of main theorem}.
}\end{exam}

\begin{rem}{\rm
A structure is said to be {\em $r$-connected} if the structure resulting from it by
removing at most $r-1$ elements (and the relationships containing them) is still connected. By modifying Theorem~{simplified form of main theorem}
and its proof, with the use of the general version of Theorems~\ref{connectivity and rigidity}
and~\ref{other asymptotic properties} from \cite{Kop}, one can prove that if a logic $L(\sigma) \geq FO(\sigma)$
can express the property
\[ \{ \mcA \in \mbS_\sigma : \mcA \text{ is $r$-connected} \} \]
where $r > 0$ is an integer, then it is not $\aleph_0$-compact.
}\end{rem}

\section{Proof of main results}\label{proof of main theorem}

\noindent
We will prove Theorem~\ref{main theorem}.
It can easily be modified to give a proof of Theorem~\ref{simplified form of main theorem}.
We assume that the signature $\sigma$ is countable and has only relation symbols.
Suppose that $Q$ and $R$ are different relation symbols and that the arity of $R$ is at least 2.
Let $\mbP \subseteq \mbS_\sigma$ be a property, so it is closed under isomorphism.
Suppose that $\mcM \in \mbS_\sigma$ is such that $Q^\mcM = \es$ and 
the following two conditions are satisfied:
\begin{itemize}
\item[(1)] If $\mcN$ is a finite $\sigma$-structure such that $\mcN$ is rigid, connected,
$Q^\mcN = \es$ and $M \cap N = \es$, then $\mcM_Q  \cup  \mcN \ \in \ \mbP$.

\item[(2)] If $\mcM' \equiv \mcM$ and 
$I$ and $J$ are sets such that $I$ is infinite and for all $i \in I$ and all $j \in J$
$\mcT_i \cong \mcT_R$, $\mcT'_j \cong \mcT'_R$ 
and the universes of any two structures from $\{\mcM'\} \cup \{\mcT_i : i \in I\} \cup \{\mcT'_j : j \in J\}$ 
are disjoint, then
\[ \mcM'_Q \ \cup \ \bigcup_{i \in I} \mcT_i \ \cup \ \bigcup_{j \in J} \mcT'_j \ \notin \ \mbP. \]
\end{itemize}

\noindent
Moreover, assume that $(L(\sigma), \models_{L(\sigma)})$ is an abstract logic such that
$L(\sigma) \geq FO(\sigma)$ and $\mbP$ can be expressed by $L(\sigma)$.
Let $\psi_\mbP \in L(\sigma)$ be a sentence  such that
for all $\mcM \in \mbS_\sigma$, $\mcM \models_{L(\sigma)} \psi_\mbP$ if and only
if $\mcM \in \mbP$. 

Theorem~\ref{main theorem} follows from Claim~A below, as will be explained.
To prove part~(ii) of Theorem~\ref{main theorem}, 
that $L(\sigma)$ does not have the $\aleph_0$-compactness property, it suffices to
find a countable set $\Gamma \subseteq L(\sigma)$ of sentences such that
\begin{itemize}
\item[(a)] if $\mcM \in \mbS_\sigma$ and $\mcM \models_{L(\sigma)} \Gamma$, then $\mcM \notin \mbP$, and
\item[(b)] for every finite $\Delta \subseteq \Gamma$ there is $\mcM \in \mbP$ such that 
$\mcM \models_{L(\sigma)} \Delta$.
\end{itemize}

\noindent
{\bf Claim A} {\em There is a countable set
$\Gamma \subseteq FO(\sigma)$ of sentences such that~(a) and~(b) hold with `$\models_{FO(\sigma)}$'
in place of `$\models_{L(\sigma)}$'.}
\\

\noindent
Observe that since $L(\sigma) \geq FO(\sigma)$, Claim~A implies part~(ii) of Theorem~\ref{main theorem}. 
Note that if $\mbP$ can be expressed by $FO(\sigma)$ and 
a countable set of sentences $\Gamma \subseteq FO(\sigma)$ satisfies~(b), then~(a) cannot
be satisfied, because of $\aleph_0$-compactness of first-order logic. 
Hence, also part~(i) of Theorem~\ref{main theorem} follows from Claim~A. 

Now suppose that $L(\sigma)$ is closed under negation.
By Claim~A and since $L(\sigma) \geq FO(\sigma)$, there is 
a countable set $\Gamma \subseteq L(\sigma)$ such that~(a) and~(b) hold.
From~(a) we get $\Gamma \models_{L(\sigma)} \neg\psi_\mbP$.
Suppose that $(\Pi, p, c)$ is a sound and complete proof system for $L(\sigma)$,
so (by completeness) there is $\pi \in \Pi$ such that $c(\pi) = \neg\psi_\mbP$,
$p(\pi) = (\varphi_1, \ldots, \varphi_k)$ and $\varphi_1, \ldots, \varphi_k \in \Gamma$.
By soundness, $\{\varphi_1, \ldots, \varphi_k\} \models_{L(\sigma)} \neg\psi_\mbP$,
which contradicts~(b).
This shows that part~(iii) of Theorem~\ref{main theorem} follows from Claim~A.

It remains to prove Claim~A.
{\em From now on we abbreviate `$\models_{FO(\sigma)}$' with `$\models$'.}
Recall the assumption that $\mcM \in \mbS_\sigma$ is such that $Q^\mcM = \es$ and~(1) and~(2) hold.
For every $\varphi \in FO(\sigma \setminus \{Q\})$ let the {\em relativization to $Q$ of $\varphi$},
denoted $\varphi_Q$, be the formula obtained from $\varphi$ by replacing every
universal quantification `$\forall x \ldots$' by `$\forall x (Q(x, \ldots, x) \rightarrow \ldots)$'
and every existential quantification `$\exists x \ldots$' by
`$\exists x (Q(x, \ldots, x) \wedge \ldots)$'.
Note that $Q$ is treated like a unary symbol although its arity may be higher than one.
Let 
\[ Th(\mcM)_Q \ = \ \big\{\varphi_Q : \varphi \in FO(\sigma \setminus \{Q\}) 
\text{ is a sentence and } \mcM \models \varphi \big\}. \]
Let $\Gamma_0$ be a set of $FO(\sigma)$-sentences  which expresses the following properties:
\begin{itemize}
\item[(3)] If $Q(x_1, \ldots, x_q)$ then $x_1 = \ldots = x_q$.

\item[(4)] If $Q(x, \ldots, x)$ and $\neg Q(y, \ldots, y)$ then $x$ and $y$ are not adjacent.

\item[(5)] The substructure whose elements $x$ are those that satisfy $\neg Q(x, \ldots, x)$ is an $R$-graph.
\end{itemize}
Let $\Gamma_1$ be the set of $FO(\sigma)$-sentences which express the following properties:
\begin{itemize}
\item[(6)] There is no cycle in the $R$-graph defined by those $x$ that satisfy $\neg Q(x, \ldots, x)$.

\item[(7)] If $\neg Q(x, \ldots, x)$ then $x$ is adjacent to exactly 4 or exactly 5 elements 
$y$ such that $\neg Q(y, \ldots, y)$.

\item[(8)] For every positive integer $k$, there are at least $k$ elements $x$ such that 
$\neg Q(x, \ldots, x)$ and $x$ is adjacent to exactly 4 elements $y$ such that $\neg Q(y, \ldots, y)$.

\item[(9)] For every positive integer $k$, if $\neg Q(x, \ldots, x)$ and $\neg Q(y, \ldots, y)$, 
$x$ is adjacent to exactly 4 elements $z$ such that $\neg Q(z, \ldots, z)$ and 
$y$ is adjacent to exactly 4 elements $u$ such that $\neg Q(u, \ldots, u)$,
then there is {\em no} $k$-path $w_1, \ldots, w_k$ such that $x = w_1$ and $y = w_k$.
\end{itemize}
Let 
\[ \Gamma \ = \ Th(\mcM)_Q \cup \Gamma_0 \cup \Gamma_1\]
and note that $\Gamma$ is countable, because $\sigma$ and hence $FO(\sigma)$, is countable. 
The definition of $\Gamma$ immediately implies the following:
\\

\noindent
{\bf Claim B} {\em If $\mcA \in \mbS_\sigma$ and $\mcA \models \Gamma$ then
\[ \mcA \ = \ \mcM'_Q \ \cup \ \bigcup_{i \in I} \mcT_i \ \cup \ \bigcup_{j \in J} \mcT'_j, \]
where $\mcM' \equiv \mcM$, $I$ is infinite, $\mcT_i \cong \mcT_R$ for every $i \in I$,
$\mcT'_j \cong \mcT'_R$ for every $j \in J$ and
the universes of any two structures from $\{\mcM'\} \cup \{\mcT_i : i \in I\} \cup \{\mcT'_j : j \in J\}$ 
are disjoint.}
\\

\noindent
It is also easy to see that if $\mcA$ has the form given by the claim,
then $\mcA \models \Gamma$, so $\Gamma$ has a model.
But the existence of a model of $\Gamma$ also follows from Claim~C below
and $\aleph_0$-compactness of $FO(\sigma)$.
From Claim~B and the assumption that $\mcM$ satisfies~(1) and~(2) it follows
that if $\mcA \models \Gamma$, then $\mcA \notin \mbP$.
So~(a) is satisfied for our choice of $\Gamma$.
It remains to prove that also~(b) is satisified. 
By assumption~(1) it suffices to prove the following:
\\

\noindent
{\bf Claim C} {\em For every finite $\Delta \subseteq \Gamma$ there is $\mcN \in \mbS_\sigma$
such that $\mcN$ is finite, rigid, connected, $N \neq \es$, $Q^\mcN = \es$,
$M \cap N = \es$ and $\mcM_Q \cup \mcN \models \Delta$.}
\\

\noindent
{\em Proof of Claim~C.}
Observe that $\mcM_Q \models Th(\mcM)_Q$ and whenever $\mcN \in \mbS_\sigma$ is an $R$-graph such that
$M \cap N = \es$, then $\mcM_Q \cup \mcN \models Th(\mcM)_Q \cup \Gamma_0$.
Therefore it suffices to find, for every finite $\Delta_1 \subseteq \Gamma_1$, an $R$-graph
$\mcN \in \mbS_\sigma$ such that $\mcN$ is finite, rigid, connected, $N \neq \es$,
$Q^\mcN = \es$ and $\mcN \models \Delta_1$.

For every integer $m \geq 3$, let $\Gamma_{1,m}$ be the subset of $\Gamma_1$ which consists
of the sentences which express the following properties:
\begin{itemize}
\item[(6')] For every $3 \leq k \leq m$, there is no $k$-cycle in the $R$-graph defined by those $x$ that satisfy 
$\neg Q(x, \ldots, x)$.

\item[(7)] If $\neg Q(x, \ldots, x)$ then $x$ is adjacent to exactly 4 or exactly 5 elements 
$y$ such that $\neg Q(y, \ldots, y)$.

\item[(8')] There are at least $m$ elements $x$ such that 
$\neg Q(x, \ldots, x)$ and $x$ is adjacent to exactly 4 elements $y$ such that $\neg Q(y, \ldots, y)$.

\item[(9')] For every $k \leq m$, if $\neg Q(x, \ldots, x)$ and $\neg Q(y, \ldots, y)$, 
$x$ is adjacent to exactly 4 elements $z$ such that $\neg Q(z, \ldots, z)$ and 
$y$ is adjacent to exactly 4 elements $u$ such that $\neg Q(u, \ldots, u)$,
then there is {\em no} $k$-path $w_1, \ldots, w_k$ such that $x = w_1$ and $y = w_k$.
\end{itemize}
Suppose that $\Delta_1 \subseteq \Gamma_1$ is finite.
Then $\Delta_1 \subseteq \Gamma_{1,m}$ for some $m$, 
so it suffices to prove that there is an $R$-graph $\mcN \in \mbS_\sigma$ such that $\mcN$ is finite, rigid, 
connected and $\mcN \models \Gamma_{1,m}$.
Because of the obvious way of transforming an undirected graph to an $R$-graph
(as described in Section~\ref{results}), it is enough to show that there exists a finite undirected graph $(V, E)$ 
with the following properties:
\begin{itemize}
\item[(6'')] For every $3 \leq k \leq m$, $(V, E)$ has no $k$-cycle.
\item[(7'')] The degree of every $v \in V$ is either 4 or 5.
\item[(8'')] At least $m$ elements (vertices) of $V$ have degree 4.
\item[(9'')] For every $k \leq m$, there does {\em not} exist a path $v_1, \ldots, v_k \in V$
such that $v_1 \neq v_k$ and both $v_1$ and $v_k$ have degree 4. 
\item[(10)] $(V, E)$ is rigid.
\item[(11)] $(V,E)$ is connected.
\end{itemize}

\noindent
For every positive integer $n$, let $\mbG_{n,5}$ be the set of 
undirected graphs $(V, E)$ without loops such that $V = \{1, \ldots, n\}$ and every
vertex $v \in V$ has degree at most $5$.
The existence of an undirected graph $(V,E)$ with the properties listed above follows from
Theorems~\ref{connectivity and rigidity} and~\ref{other asymptotic properties} 
below and the easily proved fact (which also follows from \cite{BC} or \cite{Bol80},
for example) that $\lim_{n\to\infty}|\mbG_{n,5}| = \infty$.

Corollary~2.4 in \cite{Kop}, which is a consequence of results of McKay and Wormald \cite{MW84, Wor81},
Koponen~\cite{Kop} and Theorem~4.3.4 in \cite{EF}, implies the following:

\begin{theor}\label{connectivity and rigidity}(\cite{Kop}, using \cite{MW84, Wor81})
The proportion of graphs $(V, E) \in \mbG_{n,5}$ that are connected and rigid approaches 1 as $n \to \infty$.
\end{theor}

\noindent
Theorems~2.1 and~3.1 in \cite{Kop} imply the following:

\begin{theor}\label{other asymptotic properties} \cite{Kop}
Fix an arbitrary integer $m \geq 3$.
There is a real number $\alpha_m > 0$ such that the proportion of graphs $(V, E) \in \mbG_{n,5}$ with
properties~(6'')--(9'') approaches $\alpha_m$ as $n \to \infty$.
\end{theor}

\noindent
This concludes the proof of Claim~C. Hence the proof of Claim~A is finished and therefore also 
the proof of Theorem~\ref{main theorem} is finished.

We conclude this section with some comments on the proof.
For any integer $m \geq 3$, the existence of an undirected graph $(V,E)$ that 
satisfies~(6'')--(8'') and~(10)--(11)
follows from results in \cite{Bol80, MW84, Wor81, Wor81b}.
However these results do not guarantee that such $(V,E)$ also satisfies~(9''), which
is why we use results from \cite{Kop}.
There is nothing special about considering (in~(7) and~(8)) $R$-graphs such that all vertices have degree~4 or~5.
We could as well have specified that the degrees are $d-1$ or $d$ for any $d \geq 5$ and 
considered the set $\mbG_{n,d}$ of all graphs with vertex set $V = \{1, \ldots, n\}$ such that
every vertex has degree at most $d$.

\section{Proofs of the applications without random graph theory}
\label{Proofs of the applications without random graph theory}

\noindent
Theorems~\ref{main theorem} and~\ref{simplified form of main theorem}
are convenient as they can be used in a uniform way to prove the conclusions
of all examples mentioned in Section~\ref{examples}, whether they speak
of authomorphism groups or of connected components (as the last example).
Moreover, they can be used on some combinations of the properties considered, such as
``the structure is rigid and connected''.

However, the conclusions of all examples in Section~\ref{examples} can be proved
without the use of random graph theory, as we explain in this section.
In all these applications except for Example~\ref{8th example} the conclusion
follows from a variant of Theorem~\ref{main theorem} or~\ref{simplified form of main theorem}
which is proved in a similar way as these theorems with the help of
the following basic result instead of Theorems~\ref{connectivity and rigidity} 
and~\ref{other asymptotic properties}:

\begin{lem}\label{tree construction}
Let $\Gamma$ be a first-order theory in a language with one binary relations symbol (besides
the identity symbol) which expresses the following properties:
\begin{itemize}
\item[(i)] The models are trees, i.e. undirected graphs without cycles.
\item[(ii)] For every integer $k > 0$,
\begin{itemize}
\item[(a)] the number of elements of degree 4 is at least $k$, and
\item[(b)] if the distance between $x$ and some element with degree 4 is $k$, then the degree of $x$ is 3.
\end{itemize}
\end{itemize}
Then
\begin{itemize}
\item[(1)] every finite $\Delta \subseteq \Gamma$ has a rigid model, and
\item[(2)] for every $\mcM \models \Gamma$, $Aut(\mcM)$ has a subgroup which is isomorphic with $S(\omega)$.
\end{itemize}
\end{lem}

\noindent
{\bf Proof.}
Let $\mcT$ be an infinite graph such that $\mcT$ is connected and has no cycle, $\mcT$ has exactly 
one element with degree 4 and all other elements have degree 3.
From the definition of $\Gamma$ it is clear that every $\mcM \models \Gamma$ is
a disjoint union of infinitely many structures each one of which is isomorphic with $\mcT$.
Then it is clear that $Aut(\mcM)$ has a subgroup which is isomorphic with $S(\omega)$, so we 
have proved~(2).

Let $\mcB_0$ be an undirected graph with only one element and call this element the {\em root} of $\mcB_0$
For $n \geq 0$ let $B_{n+1}$ be the graph constructed in the following way:
Take two disjoint copies $\mcB'_n$ and $\mcB''_n$ of $\mcB_n$ and a new element $b_{n+1}$ and
let $\mcB_{n+1}$ be the tree obtained by adding an edge between $b_{n+1}$ and the root of $\mcB'_n$
as well as an edge between $b_{n+1}$ and the root of $\mcB''_n$; call $b_{n+1}$ the root of $\mcB_{n+1}$.
Hence, $\mcB_n$ is a binary tree of height $n$ for each $n \geq 0$.
The elements of a tree with degree 1 are called {\em leaves}.

For each integer $k > 0$ let $\mcP_k$ denote a path of length $k$, that is,
$P_k = \{a_0, \ldots, a_k\}$ where $a_i \sim_\mcP a_j$ if and only if $i+1 = j$.
We also assume that $P_k$ and $P_n$ are disjoint if $k \neq n$.
For each $n > 0$ let $\mcF_n$ be a forest (i.e. disjoint union of trees) constructed as follows. 
Take four disjoint copies of $\mcB_n$ and join their roots (by adding four edges) with
a new vertex and call the new tree $\mcT_n$. 
Note that $\mcT_n$ has $4 \cdot 2^n = 2^{n+2}$ leaves.
Now take disjoint copies $\mcT_{n,1}, \ldots, \mcT_{n,n}$ of $\mcT_n$ and
let $\mcF'_n$ be the union of these.
Enumerate all leaves in $\mcF'_n$ as $a_1, \ldots, a_{n2^{n+2}}$.
For each $k = 1, \ldots, n2^{n+2}$ let $\mcP'_k$ be a copy of $\mcP_k$ such
that $P'_k$ is disjoint from $P'_m$ if $k \neq m$.
For $k = 1, \ldots, n2^{n+2}$, add an edge between $a_k$ and
one end of $\mcP'_k$. Call the resulting forest $\mcF_n$.
Then it is easy to see that $\mcF_n$ is rigid.
Morover, if $\Delta \subseteq \Gamma$ is finite then $\mcF_n \models \Delta$ if
$n$ is choosen sufficiently large, so~(1) is proved.
\hfill $\square$
\\

\noindent
The conclusion of Example~\ref{8th example} can be proved by using the following 
modification of Lemma~\ref{modification of tree construction}
instead of Theorems~\ref{main theorem} and~\ref{simplified form of main theorem}:

\begin{lem}\label{modification of tree construction}
Let $\Gamma$ be a first-order theory in a language with one binary relations symbol (besides
the identity symbol) which expresses the following properties:
\begin{itemize}
\item[(i)] The models are trees.
\item[(ii)] For every integer $k > 0$,
\begin{itemize}
\item[(a)] the number of elements of degree 4 is at least $k$, and
\item[(b)] if the distance between $x$ and some element with degree 4 is $k$, then the degree of $x$ is 2.
\end{itemize}
\end{itemize}
Then
\begin{itemize}
\item[(1)] every finite $\Delta \subseteq \Gamma$ has a connected model, and
\item[(2)] every $\mcM \models \Gamma$ has infinitely many components.
\end{itemize}
\end{lem}

\noindent
{\bf Proof.} Modification of the proof of Lemma~\ref{tree construction}. 
\hfill $\square$

\section{The property:
``the automorphism group has cardinality $\leq 2^{\aleph_0}$''}

\begin{theor}\label{a compact logic}
Let $\sigma$ be a countable signature. 
Extend $FO(\sigma)$ with a new atomic sentence $P$.
Let $L(\sigma)$ be the least set of formulas such that it contains
all atomic first-order formulas and $P$ and is closed under the first-order operations.
For every $\mcM \in \mbS_\sigma$, 
let $\mcM \models_{L(\sigma)} P$ if and only if $|Aut(\mcM)| \leq 2^{\aleph_0}$.
Otherwise `$\models_{L(\sigma)}$' is defined just as for first order formulas.
Then $L(\sigma) > FO(\sigma)$ and $L(\sigma)$ has the $\aleph_0$-compactness property.
\end{theor}

\noindent
{\bf Proof.}
Clearly, $L(\sigma) \geq FO(\sigma)$.
Suppose that $\varphi \in FO(\sigma)$ is such that for all $\mcM \in \mbS_\sigma$,
$\mcM \models_{L(\sigma)} P$ if and only if $\mcM \models_{L(\sigma)} \varphi$.
By the use of Ehrenfeucht-Mostowski models, in particular we can use Lemma~5.2.7
in \cite{Mar}, there is $\mcM \models_{L(\sigma)} \neg P$ and hence 
$\mcM \models_{L(\sigma)} \neg\varphi$, so $\mcM \models_{FO(\sigma)} \neg\varphi$.
As $FO(\sigma)$ has the $\aleph_0$-Löwenheim-Skolem property, there is a countable
model $\mcN \models_{FO(\sigma)} \neg\varphi$. 
Then also $\mcN \models_{L(\sigma)} \neg\varphi$ so $\mcN \models_{L(\sigma)} \neg P$,
which is impossible because every countable structure has at most $2^{\aleph_0}$ automorphisms.
Hence $L(\sigma) > FO(\sigma)$. 
It remains to prove that $L(\sigma)$ has the $\aleph_0$-compactness property.

Suppose that $\Gamma \subseteq L(\sigma)$ has the property that every finite
subset of $\Gamma$ has a model. We will show that also $\Gamma$ has a model.
For every $\varphi \in L(\sigma)$, we define $\varphi' \in FO(\sigma)$-formula as follows:
\begin{itemize}
\item[(i)] If $\varphi$ is an atomic first-order formula, then $\varphi' = \varphi$.
\item[(ii)] If $\varphi = P$, then $\varphi' = \forall x (x = x)$.
\item[(iii)] If $\varphi = \neg\psi$, then $\varphi' = \neg\psi'$.
\item[(iv)] If $\varphi = \psi \wedge \theta$, then $\varphi' = \psi' \wedge \theta'$
(and similarly for the connectives $\vee, \rightarrow, \leftrightarrow$ if we use them).
\item[(v)] If $\varphi = \exists x \psi$, then $\varphi' = \exists x \psi'$ (and similarly
for $\forall$ if we use it).
\end{itemize}
For every $\varphi \in L(\sigma)$ we define $\varphi'' \in FO(\sigma)$ in a similar way,
by keeping~(i) and the inductive steps~(iii)--(v), but replacing~(ii) with
\begin{itemize}
\item[(ii')] If $\varphi = P$, then $\varphi'' = \neg\forall x (x = x)$.
\end{itemize}

\noindent
{\bf Claim.} {\em 
(a) If $|Aut(\mcM)| \leq 2^{\aleph_0}$, then for every sentence $\varphi \in L(\sigma)$,
$\mcM \models_{L(\sigma)} \varphi$ if and only if $\mcM \models_{FO(\sigma)} \varphi'$.\\
(b) If $|Aut(\mcM)| > 2^{\aleph_0}$, then for every sentence $\varphi \in L(\sigma)$,
$\mcM \models_{L(\sigma)} \varphi$ if and only if $\mcM \models_{FO(\sigma)} \varphi''$.}

\medskip
\noindent
{\bf Proof.} Easy induction. 
\hfill $\square$
\\

\noindent
Let $\Gamma' = \{\varphi' : \varphi \in \Gamma\}$ and
$\Gamma'' = \{\varphi'' : \varphi \in \Gamma\}$.
Now we are ready to show that $\Gamma$ has a model.
There are three cases to consider.

{\em Case 1: Neither $\Gamma'$ nor $\Gamma''$ has an infinite model.}
We first show that then $\Gamma'$ has a model.
For a contradiction, suppose that $\Gamma'$ does not have a model.
By the $\aleph_0$-compactness property of $FO(\sigma)$,
there is a finite $X \subseteq \Gamma$ such that $X' = \{\varphi' : \varphi \in X\}$
does not have a model.
But then, by~(a) of the claim, for every finite $Y$ such that $X \subseteq Y \subseteq \Gamma$,
$Y$ does not have a model $\mcM$ such that $|Aut(\mcM)| \leq 2^{\aleph_0}$.
Then, by the assumption that every finite subset of $\Gamma$ has a model,
for every finite $Y$ such that $X \subseteq Y \subseteq \Gamma$,
$Y$ has a model $\mcM_Y$ such that $|Aut(\mcM)| > 2^{\aleph_0}$.
In particular, $\mcM_Y$ is infinite.
Thus, by~(b) of the claim, for every finite $Y$ such that $X \subseteq Y \subseteq \Gamma$,
$\mcM_Y$ is an infinite model of $Y'' = \{\varphi'' : \varphi \in Y\}$.
By $\aleph_0$-compactness of $FO(\sigma)$, $\Gamma''$ has an infinite model,
contradicting the assumption we made.
Hence we conclude that $\Gamma'$ has a model $\mcM$.
By our assumption $\mcM$ is finite, so $|Aut(\mcM)| \leq 2^{\aleph_0}$.
Therefore part~(a) of the claim gives $\mcM \models_{L(\sigma)} \Gamma$.

{\em Case 2: $\Gamma'$ has an infinite model.}
By the $\aleph_0$-Löwenheim-Skolem property of $FO(\sigma)$,
$\Gamma'$ has a countable model $\mcM$.
Since $|Aut(\mcM)| \leq 2^{\aleph_0}$, it follows from part~(a) of the claim that
$\mcM \models \Gamma$.

{\em Case 3: $\Gamma''$ has an infinite model.}
By the use of Ehrenfeucht-Mostowski models, for example Corollary~5.2.7 in \cite{Mar},
it follows that $\Gamma''$ has a model $\mcM$ such that $|Aut(\mcM)| > 2^{\aleph_0}$.
Now part~(b) of the claim implies that $\mcM$ is a model of $\Gamma$.
\hfill $\square$

\end{document}